\documentclass[12pt]{amsart}
\usepackage[a4paper,left=3.5cm,right=3.5cm,top=4cm,bottom=4cm]{geometry}

\usepackage{amsmath,amssymb,latexsym, graphicx, mathrsfs, manfnt, textcomp, appendix, chemarrow, enumerate}
\usepackage[misc,geometry]{ifsym}
\usepackage{hyperref}
\usepackage[all]{xy}
\usepackage[utf8]{inputenc}
\usepackage[T1]{fontenc}

\numberwithin{equation}{section}

\newcommand\ptensor{\widehat{\otimes}_{\pi}}
\newcommand\itensor{\check{\otimes}_{\varepsilon}}

\newcommand\nbd{\mathscr{U}}

\newcommand{\lcp}{\mathfrak{s}}

\theoremstyle{thmit} 
\newtheorem{theorem}{Theorem}
\newtheorem{lemma}{Lemma}
\newtheorem{corollary}{Corollary}
\newtheorem{proposition}{Proposition}

\newtheorem{example}{Example}
\begin{document}
\author{Ersin Kızgut}

\address{Çankaya University, Institute of Natural and Applied Sciences, Öğretmenler Cd No.14 06530 Ankara Turkey}

\email{ersinkizgut@cankaya.edu.tr}

\author{Murat Yurdakul}

\address{Middle East Technical University, Department of Mathematics, Dumlupınar Blv No. 1 06800 Ankara Turkey}

\email{myur@metu.edu.tr}

\title{Remarks on strictly singular operators}


\maketitle

\begin{abstract}
	A continuous linear operator $T:E \to F$ is called strictly singular if it cannot be invertible on any infinite dimensional closed subspace of its domain. In this note we discuss sufficient conditions and consequences of the phenomenon $LB(E,F)=L_s(E,F)$, which means that every continuous linear bounded operator defined on $E$ into $F$ is strictly singular. 
\end{abstract}

\section{Introduction and first results}

	A continuous linear operator mapping a locally convex space (lcs) $E$ into a lcs $F$ is said to be bounded if it maps a neighborhood of the origin of $E$ into a bounded subset of $F$. If $E$ or $F$ is a normed space, then every continuous linear operator between $E$ and $F$ is bounded. A continuous linear operator $T:E \to F$ is called strictly singular if it fails to be invertible on any infinite dimensional closed subspace of $E$. On the class of Banach spaces, Kato \cite{Kat58} introduced strictly singular operators in connection with perturbation theory of Fredholm operators. A compact operator is strictly singular, while the converse is not true in general. In terms of lcs's, van Dulst first \cite{Dul70}, \cite{Dul71} studied strictly singular operators on Ptak (or $B$-complete) spaces and generalized Hilbert spaces. Wrobel \cite{Wro78} characterized strictly singular operators on lcs's for the class of $B_r$-complete spaces. $L_s(X,Y)$ is an operator ideal, if the pair $(X,Y)$ belongs to the class of Banach spaces. As proved in \cite{Sus81}, this is not the case when it belongs to the class of general lcs's. However, each of the classes of bounded strictly singular operators $LB_s(E,F)$ and of compact operators $L_c(E,F)$ forms an operator ideal in lcs's. By \cite{Ona91}, for Fréchet spaces $E$ and $F$, $LB_s(E,F)=L_s(E,F)$ if $F$ has the property $(y)$. If $E$ contains $\ell^1$ as a quotient, then $LB_s(E,F) \neq L_s(E,F)$.

\begin{lemma}\label{ss ideal in br}
	Let $E$ and $F$ be lcs's where $E$ is  $B_r$-complete. Then, $L_s(E,F)$ forms an operator ideal.	
\end{lemma}

\begin{proof}
	Suppose that $T: E \to F$ and $S:E \to F$ are strictly singular operators. Then, for any $M \leq E$, by \cite[Theorem 1-IV]{Wro78}, find $N \leq M$ such that $T|_N$ is precompact. Then find $P \leq N$ such that $S|_P$ is precompact. The ideal property of precompact operators on lcs's yields the result.
\end{proof}

We give the following proposition as an application of operator ideal property of strictly singular operators on $B_r$-complete lcs's. It is a generalization of \cite[Problem 4.5.2]{Abr02}, and in particular, it is also true when re-stated for bounded strictly singular operators acting on general lcs's.

\begin{proposition}\label{direct sum ss}
Let $E:=\bigoplus^n E_i$ and $F:=\bigoplus^m F_j$ be lcs's where $E$ is  $B_r$-complete. Then, $T: E \to F$ is strictly singular iff each of $T_{ij}:E_i \to F_j$ is strictly singular for each $i=1,2,\dots,n$ and for each $j=1,2,\dots,m$.
\end{proposition}

\begin{proof}
Assume that each $T_{ji}$ is strictly singular. Let $\pi_i:E \to E_i$ be the canonical projection and define $\rho_j:F_j \to F$ by $\rho y_j=0 \oplus 0 \oplus \dots \oplus 0 \oplus y_j \oplus 0 \oplus \dots \oplus 0$, for which $y_j$ is the $j$-th summand. Consider $E \xrightarrow{\pi_i} E_i \xrightarrow{T_{ji}} F_j \xrightarrow{\rho_j} F$, and write $S_{ji}:=\rho_j \circ T_{ji} \circ \pi_i$. Then $S_{ji}(x_1 \oplus x_2 \oplus \dots \oplus x_n)=0 \oplus 0 \oplus \dots \oplus 0 \oplus T_{ji}x_i \oplus 0 \oplus \dots \oplus 0$, where $T_{ji}$ is the $j$-th summand. By Lemma~\ref{ss ideal in br} and rewriting $T=\sum_{i=1}^n \sum_{j=1}^m S_{ji}$, $T$ is strictly singular. For the converse, let $T \in L_s(E,F)$, and suppose that the operator $T_{ji}$ is not strictly singular for some $i,j$, and for $M \leq E$, $r \in I$ and $s \in J$, $N_{rs}(T|_M):=\sup\{q_s(Tx):p_r(x) \leq 1, x \in M\}$. Then by \cite[Theorem 2.1]{Moo10}, for any $M \leq E_i$ and for some $s \in J$, $N_{rs}(T|_N) > \varepsilon$, for all $r \in I$. If we write $\widehat{M}:=\{0\} \oplus \{0\} \oplus \dots \oplus M \oplus \{0\} \oplus \dots \oplus \{0\}$ where $M$ places in the $i$-th summand, $\widehat{M}$ is a vector subspace of $X$. $N_{rs}(\tau_{\widehat{M}})>\varepsilon$, for all $r \in I$. But, that contradicts the assumption $T$ is strictly singular. 
\end{proof}

Every compact operator is clearly bounded. If we replace compactness with strict singularity, the results in \cite{Ona91}, \cite{Yur93} and \cite{Ter86} show that a strictly singular operator $T:E \to F$ between lcs's cannot be unbounded, since such an assumption leads to a contradiction with the well-known result of Bessaga, Pelczynski, and Rolewicz \cite{Bes61}. In Section~\ref{section: sufficient_conditions_ss} we survey the additional hypotheses under which a bounded operator between two lcs's is automatically strictly singular.

\section{Strict singularity of bounded operators}\label{section: sufficient_conditions_ss}

In this part, the additional hypotheses for strict singularity of a bounded operator will be introduced. Our starting point will be the class of Banach spaces. Applications of \cite[Lemma 2]{Dja98-2} shall yield that some of these results are meaningful on the class of lcs's. In a particular case, to obtain a characterization (Theorem~\ref{characterization}) is possible. The most commonly known non-trivial example for $L(X,Y) = L_s(X,Y)$ in Banach spaces is when $X=\ell^p$ and $Y=\ell^q$ such that $1 \leq p <q<\infty$. This phenomenon was essential in the isomorphic classification of Cartesian products of power series spaces \cite{Dja98-2}, and in the problem whether the sum of two complemented subspaces is also complemented \cite{Kiz15} (if $X$ and $Y$ be complemented subspaces of $E$ then, $X+Y$ is a complemented subspace of $E$ if $LB(X,Y)=L_s(X,Y)$). A Banach space $X$ is said to possess the Schur property (SP) if every weakly convergent sequence in $X$ converges in norm. $X$ is called almost reflexive, if every bounded sequence is weakly Cauchy. $X$ is called weakly sequentially complete (wsc) if every weakly Cauchy sequence in $X$ is weakly convergent. $X$ is said to have the Dieudonné property (DP) if operators on $X$ mapping weakly Cauchy sequences into weakly convergent sequences are weakly compact.

A pair of Banach spaces $(X,Y)$ is called totally incomparable if there exists no Banach space $Z$ which is isomorphic to a subspace of $X$ and to a subspace of $Y$. A property $P$ on a Banach space $X$ is called hereditary if it is enjoyed by every $M \leq X$. $X$ is said to have nowhere $P$ if it has no subspace having the property $P$. By $L_w$ and $L_v$ we denote the classes of weakly compact and completely continuous (or fully complete) operators, respectively. A Banach space $X$ is said to have the Dunford-Pettis property (DPP) if $L_w(X,Y) \subseteq L_v(X,Y)$, for any Banach space $Y$. $X$ is said to have the reciprocal Dunford-Pettis property (rDPP), if $L_v(X,Y) \subseteq L_w(X,Y)$. 

\begin{lemma}\label{lemma: aux}
Let $X$ and $Y$ be Banach spaces.
	\begin{enumerate}[1.]
		\item Let $X$ be almost reflexive. Then, for any wsc Banach space $Y$, $L(X,Y) = L_s(X,Y)$.
		\item Let $X$ have DP, and let $Y$ be wsc. Then, $L(X,Y) = L_w(X,Y)$.
		\item Let $Y$ be almost reflexive. Then $L(X,Y) = L_v(X,Y)$ implies $L(X,Y) = L_s(X,Y)$.
		\item Let $X$ be a Banach space with SP. Then, for every $M \leq X$, $\ell^1 \hookrightarrow M$.
		\end{enumerate}
\end{lemma}

\begin{proof}
	\begin{enumerate}[1.]
		\item Since $X$ is almost reflexive, if $(x_n)$ is a bounded sequence in $X$, then $(Tx_n)$ has a weakly Cauchy sequence in $Y$. But $Y$ is wsc, that is, every weakly Cauchy sequence converges weakly in $Y$. Therefore, $T$ is weakly compact. 
		\item Let $(x_n) \in X$ be weakly Cauchy, and let $T \in L(X,Y)$. Then $(Tx_n)$ is weakly Cauchy in $Y$. Since $Y$ is assumed to be wsc, $(Tx_n)$ converges weakly. But $X$ has DP, so $T \in L_w(X,Y)$.
		\item See \cite[Theorem 1.7]{Lac65}.
		\item Let $X$ have the SP and suppose there exists $M \leq X$ not containing $\ell^1$. Then, any bounded sequence $(x_n)$ in $M$, has a weakly Cauchy subsequence since $M$ is equivalently almost reflexive. However, $M$ inherits SP. Then the weakly Cauchy subsequence of $(x_k)$ converges in $X$. Therefore, $M$ is finite dimensional. Contradiction. 
	\end{enumerate}
\end{proof}

\begin{theorem}\label{theorem: sufficient conditions}
				Let $X$ and $Y$ be Banach spaces. Each of the following implies $L(X,Y) = L_s(X,Y)$.
					\begin{enumerate}[1.]
						\item $X$ and $Y$ are totally incomparable.
						\item $X$ is nowhere reflexive, $Y$ is reflexive.
						\item $X$ is nowhere reflexive, $Y$ is quasi-reflexive. 
						\item $X$ is almost reflexive and nowhere reflexive, $Y$ is wsc (see Example~\ref{vir_ir}).
						\item $Y$ has hereditary $P$, $X$ has nowhere $P$.
						\item $Y$ is almost reflexive, $X$ is hereditarily-$\ell^1$.
						\item $X$ has SP, $Y$ is almost reflexive.
						\item $X$ is reflexive, $Y$ has SP.
						\item $X$ has the hereditary DPP, $Y$ is reflexive.
						\item $L(X,Y) = L_w(X,Y)$ and $X$ has DPP.
						\item $L(X,Y) = L_v(X,Y)$ and $X$ has rDPP.
						\item $X$ is a Grothendieck space with DPP, $Y$ is separable.
						\item $X$ has both DP and DPP, $Y$ is wsc (see Example~\ref{vir_ir}).
					\end{enumerate}
			\end{theorem}
			
\begin{proof}
	\begin{enumerate}[1.]
		\item Suppose $T:X \to Y$ is a non-strictly singular operator. So find  $M \leq X$ on which $M \simeq T(M) \leq Y$. Since $X$ and $Y$ are totally incomparable, this is impossible.
		\item See \cite[Theorem b]{Gol63}. This result is generalized in part 5.
		\item Suppose there exists a non-strictly singular operator $T \in L(X,Y)$. Then, $T$ is an isomorphism when restricted to  $M \leq X$, so $M$ is quasi-reflexive. However, by \cite[Lemma 2]{Her67} there exists a reflexive $N \leq M$. This contradicts the assumption $X$ is nowhere reflexive.
		\item By part 1 of Lemma~\ref{lemma: aux}, $L(X,Y) = L_w(X,Y)$. Now let $T: X \to Y$ which has a bounded inverse on $M \leq X$. If $(x_n)$ is a bounded sequence in $M$, then there exists $(T{x_k}_n)$ a weakly convergent subsequence of $(Tx_n)$ in $Y$. Hence $({x_k}_n)$ is weakly convergent in $M$, since $T$ has a bounded inverse on $M$. Thus, every bounded sequence in $M$ has a weakly convergent subsequence in $M$. So $M$ is reflexive. Contradiction.
		\item For some $M \leq X$ suppose there exists $T:X \to Y$ such that $M \simeq T(M)$. But $T(M)$ inherits $P$. Hence $M$ has $P$. This contradicts $X$ has nowhere $P$. Now let $S: Y \to X$ be such that $N \simeq S(N) \subseteq X$ for some $N \leq Y$. Since $X$ has nowhere $P$, $S(N)$ does not enjoy $P$. Contradiction.
		\item A particular case of part 5.
		\item Any operator $T$ with range $Y$ maps bounded sequences into weakly Cauchy sequences, since $Y$ is almost reflexive. On the other hand, any such $T$ defined on $X$ maps weakly Cauchy sequences into norm convergent sequences by the SP. That implies $L(X,Y) = L_v(X,Y)$. Hence by part 3 of Lemma~\ref{lemma: aux} the result follows.
		\item Let $T \in L(X,Y)$ have a bounded inverse on some $M \leq X$, that is, $M \simeq T(M)$. Since $Y$ have SP, it also has the hereditary DPP \cite{Die80}. Hence so does $M$. But $M$ is reflexive. By \cite[Theorem 2.1]{Kes72}, reflexive spaces have nowhere DPP. Contradiction.
		\item Let $M \leq X$ on which an arbitrary operator $T:X \to Y$ has a bounded inverse. Then, $M \simeq T(M)$. So, $M$ is reflexive. Hence $M$ cannot have DPP. Contradiction.
		\item Since $X$ has DPP, $L(X,Y) \in L_v(X,Y)$. Then, by \cite[Theorem 2.3]{Lac65}, the result follows.
		\item Since $X$ has the rDPP and $L(X,Y) = L_v(X,Y)$, $T \in L_v \cap L_w$. By \cite[Theorem 2.3]{Lac65}, we are done.
		\item By \cite[Theorem 4.9]{Mor00}, any such operator $T:X \to Y$ is weakly compact. Since $X$ has the DPP, $T$ is completely continuous. By part 10, we reach the result.
		\item By part 2 of Lemma~\ref{lemma: aux}, $L(X,Y) = L_w(X,Y)$. $X$ possesses the DPP, so $L(X,Y) = L_v(X,Y)$. By \cite[Theorem 2.3]{Lac65}, the proof is completed.
	\end{enumerate}
\end{proof}

\begin{example}\label{vir_ir}
\textup{Note that the non-reflexive space $c_0$ is almost reflexive. Suppose there exists a reflexive subspace $E$ of $c_0$. Since $c_0$ fails SP, it is not isomorphic to any subspace of $E$. But this contradicts \cite[Proposition 2.a.2]{Lin96}. The space $C(K)$, where $K$ is a compact Hausdorff space enjoys both DP \cite{Gro53}, and DPP \cite{Kes72}}.
\end{example}

\begin{corollary}
				Let $X,Y,W,Z$ be Banach spaces. Then,
				\begin{enumerate}[1.]
					\item If $X',Y',Z'$ have SP and $W$ is almost reflexive, every operator defined on $(X \ptensor Y)'$ into $W \ptensor Z$ is strictly singular.
					\item If $X$ and $Y$ are reflexive spaces one of which having the approximation property, $L(X,Y') = L_c(X,Y')$, $W$ and $Z$ have SP, then every operator defined on $X \ptensor Y$ into $W \itensor Z$ is strictly singular.
					\item If $X$ is almost reflexive and $Y$ is almost reflexive and has DPP, then every operator defined on $X \ptensor Y$ into $\ell^1$ is strictly singular.
				\end{enumerate}
			\end{corollary}

\begin{proof}
	\begin{enumerate}[1.]
		\item 	By \cite[Corollary 1.6]{Lac65}, $L(W,Z') = L_c(W,Z')$. So by \cite[Theorem 3]{Emm92} we deduce $W \ptensor Z$ is almost reflexive. On the other hand, by \cite[Theorem 3.3(b)]{Rya87} we reach that $L(X,Y')$ has SP. But in \cite{Rya02} it is proved that $L(X,Y') \simeq (X \ptensor Y)'$. So $(X \ptensor Y)'$ has SP. Therefore, Thorem~\ref{theorem: sufficient conditions} part 7 yields the result.
		\item 	By \cite[Theorem 4.21]{Rya02}, $X \ptensor Y$ is reflexive. By \cite{Lus75}, SP respects injective tensor products. So $W \itensor Z$ has SP. Then, part 8 of Theorem~\ref{theorem: sufficient conditions} finishes the proof.
		\item By \cite{Die80}, $Y'$ has SP. Then by \cite[Corollary 1.6]{Lac65}, $L(X,Y') = L_c(X,Y')$. Hence, \cite[Theorem 3]{Emm92} yields that $X \ptensor Y$ is almost reflexive. It is clear that every operator defined from an almost reflexive space into $\ell^1$ is strictly singular.
	\end{enumerate}
\end{proof}
			
Now let us turn our attention to general class of lcs's. Let $\mathsf{P}$ be a class of Banach spaces having a certain hereditary property $P$. Then, $\lcp(\mathsf{P})$ \cite{Cas91} is defined by the set of lcs's $E$ with local Banach spaces $E_U \in \mathsf{P}$ for which $E_U$ is the completion of the normed space obtained by $E/{p_u^{-1}(0)}$, where $U \in \nbd(E)$ and $p_u$ its gauge functional. For instance, by $\lcp(\mathsf{X})$, we denote the class of lcs's such that each of their local Banach spaces are hereditarily $\ell^1$. A lcs $E$ is called locally Rosenthal \cite{Boy07}, if it can be written as a projective limit of Banach spaces each of which contains no isomorpic copy of $\ell^1$. A quasinormable lcs $E$ admitting no isomorphic copies of $\ell^1$ is locally Rosenthal. By $\lcp(\mathsf{V})$, we denote the class of lcs's with local Banach spaces each of which having SP. A lcs $E$ is called infra-Schwartz (or Komura) if each of its local Banach spaces is reflexive. An infra-Schwartz space turns out to be locally Rosenthal, as proved in \cite{Boy07}. In \cite{Dja98-2}, it is shown that the relation $L(X_k,Y_m)=L_s(X_k,Y_m)$ for every $k,m$ is sufficient for $L(X,Y)=L_s(X,Y)$, where $X= \projlim_k X_k$ and $Y=\projlim_m Y_m$.

Let $\lambda_1(A) \in (d_2)$, and $\lambda_p(A) \in (d_1)$ as in \cite{Dra65}. Then, by \cite{Zah73}, $L(\lambda_1(A),\lambda_p(A)) = LB(\lambda_1(A),\lambda_p(A))$. For $1 \leq p < \infty$, we know that $\displaystyle \lambda_p(A)=\projlim \ell^p(a_n)$. Since $\ell^p(a_n), 1<p<\infty$ has no subspace isomorphic to $\ell^1$, $L(\ell^1,\ell^p) = L_s(\ell^1,\ell^p)$. Then, by \cite[Lemma 2]{Dja98-2}, $L(\lambda_1(A),\lambda_p(A)) = L_s(\lambda_1(A),\lambda_p(A))$. Resting on the same argument, to obtain several sufficient conditions for $L(E,F)=L_s(E,F)$ is possible for the class of general lcs's.

\begin{theorem}\label{locRos}
Let $E,F$ be lcs's. Each of the following implies $LB(E,F)=L_s(E,F)$.

\begin{enumerate}[1.]
	\item $E \in \lcp(\mathsf{X})$ and $F$ is locally Rosenthal.
	\item $E \in \lcp(\mathsf{V})$ and F is a quasinormable Fréchet space.
	\item $E$ is infra-Schwartz and $F \in \lcp(\mathsf{V})$
	\item $E \in \lcp(\mathsf{P}^\neg)$ and $F \in \lcp(\mathsf{P})$
\end{enumerate}
\end{theorem}

\begin{proof}
\begin{enumerate}[1.]
\item Since $F$ is locally Rosenthal, there exists a family of Banach spaces $\{F_m\}$ each of which does not contain an isomorphic copy of $\ell^1$ such that $F=\projlim F_m$. Because $E \in \lcp(\mathsf{X})$, there exists a family of Banach spaces $\{E_k\}$ such that every $M_k \leq E_k$ contains a subspace isomorphic to $\ell^1$. By part 6 of Theorem~\ref{theorem: sufficient conditions}, any linear operator $T_{mk}: E_k \to F_m$ is strictly singular. Making use of \cite[Lemma 2]{Dja98-2}, we reach the result.
\item By \cite[Theorem 6]{Min95}, $F$ is locally Rosenthal. Since $E \in \lcp(\mathsf{V})$, by part 4 of Lemma~\ref{lemma: aux}, $E \in \lcp(\mathsf{X})$. Then, by part 1, we are done.
\item Since $E$ is infra-Schwartz, any of its local Banach spaces $E_k$ is reflexive. The assumption on $F$ completes the conditions in part 8 of Theorem~\ref{theorem: sufficient conditions}. Combined with \cite[Lemma 2]{Dja98-2}, we are done.
\item Since $E \in \lcp(\mathsf{P}^\neg)$, one may rewrite $E= \projlim_k E_k$, where each $E_k$ has no subspace having property $P$. Similarly, $F=\projlim_m F_m$ where each $F_m$ is hereditarily $P$. Hence, by part 5 of Theorem~\ref{theorem: sufficient conditions}, $L(E_k,F_m) = L_s(E_k,F_m)$ for every $k,m$. Applying \cite[Lemma 2]{Dja98-2}, we obtain $LB(E,F) = L_s(E,F)$.
\end{enumerate}
\end{proof}

\begin{theorem}\label{characterization}
	Let $(E,F,G)$ be a triple of Fréchet spaces satisfying the following
\begin{enumerate}[1.]
\item Every subspace of $E$ contains a subspace isomorphic to $G$.
\item $F$ has no subspace isomorphic to $G$.
\end{enumerate}
Then, $LB(E,F) = L_s(E,F)$. Let $F$ have continuous norm in addition. Then, (2) is also necessary if $F$ is a Fréchet-Montel space.
\end{theorem}

\begin{proof}
The sufficiency part is very similar to the proof of Theorem~\ref{locRos} part 4. For necessity, let $E$ be a Fréchet space and let $F$ be an (FM)-space admitting a continuous norm. Let any linear operator $T:E \to F$ be strictly singular. Then, by \cite[Proposition 1]{Yur93}, it is bounded. Now let there exist $N \leq Y$ which is isomorphic to $G$. Then $I|_N:N \to G$ is bounded, hence compact. Then $N$ is finite dimensional. Contradiction. 
\end{proof}

\bibliographystyle{plain}
\bibliography{ss}

\begin{thebibliography}{10}

\bibitem{Abr02}
Y.~A. Abramovich and C.~D. Aliprantis.
\newblock {\em Problems in Operator Theory}.
\newblock Graduate Studies in Mathematics. American Mathematical Society, 2002.

\bibitem{Bes61}
C.~Bessaga, A.~Pelczynski, and S.~Rolewicz.
\newblock On diametral approximative dimension and linear homogeneity of
  {F}-spaces.
\newblock {\em Bull. Acad. Polon. Sci}, 9:677--683, 1961.

\bibitem{Boy07}
C.~Boyd and M.~Venkova.
\newblock Grothendieck space ideals and weak continuity of polynomials on
  locally convex spaces.
\newblock {\em Monatsch. Math}, 151:189--200, 2007.

\bibitem{Cas91}
J.~M.~F. Castillo and M.~A. Sim{\~o}es.
\newblock Some problems for suggested thinking in {F}r{\'e}chet space theory.
\newblock {\em Extr. Math}, 6:96--114, 1991.

\bibitem{Sus81}
S.~Dierolf.
\newblock A note on strictly singular and strictly cosingular operators.
\newblock {\em Indag. Math}, 84:67--69, 1981.

\bibitem{Die80}
J.~Diestel.
\newblock A survey of results related to the {D}unford-{P}ettis property.
\newblock {\em Contemp. Math}, 2:15--60, 1980.

\bibitem{Dja98-2}
P.~B. Djakov, S.~{\"O}nal, T.~Terzio{\~g}lu, and M.~Yurdakul.
\newblock Strictly singular operators and isomorphisms of {C}artesian products
  of power series spaces.
\newblock {\em Arch. Math}, 70:57--65, 1998.

\bibitem{Dra65}
M.~M. Dragilev.
\newblock On regular basis in nuclear spaces.
\newblock {\em Math. Sbornik}, 68:153--175, 1965.

\bibitem{Emm92}
G.~Emmanuelle.
\newblock Banach spaces in which {D}unford-{P}ettis sets are relatively
  compact.
\newblock {\em Arch. Math}, 58:477--485, 1992.

\bibitem{Gol63}
S.~Goldberg and E.~O. Thorp.
\newblock On some open questions concerning strictly singular operators.
\newblock {\em Proc. Amer. Math. Soc}, 14:224--226, 1963.

\bibitem{Gro53}
A.~Grothendieck.
\newblock Sur les applications lineares faiblement compactes d'espaces du type
  {C}({K}).
\newblock {\em Canad. J. Math}, 5:129--173, 1953.

\bibitem{Her67}
R.~Herman and R.~J. Whitley.
\newblock An example concerning reflexivity.
\newblock {\em Studia Math}, 28:289--294, 1967.

\bibitem{Kat58}
T.~Kato.
\newblock Perturbation theory for nullity, deficiency and other quantities of
  linear operators.
\newblock {\em J. Analyse Math}, 6:261--322, 1958.

\bibitem{Kes72}
M.~C. Kester.
\newblock {\em The {D}unford-{P}ettis property}.
\newblock PhD thesis, Oklahoma State University, 1972.

\bibitem{Kiz15}
E.~K{\i}zgut and M.~Yurdakul.
\newblock Some stability results in projective tensor products.
\newblock arxiv: 1504.03254, 2015.

\bibitem{Lac65}
E.~Lacey and R.~J. Whitley.
\newblock Conditions under which all the bounded linear maps are compact.
\newblock {\em Math. Ann}, 58:1--5, 1965.

\bibitem{Lin96}
J.~Lindenstrauss and L.~Tzafriri.
\newblock {\em Classical Banach Spaces I and II}.
\newblock Springer-Verlag, 1996.

\bibitem{Lus75}
F.~Lust.
\newblock Produits tensoriels injectifs d'espaces de sidon.
\newblock {\em Colloq. Math}, 32:286--289, 1975.

\bibitem{Min95}
M.~A. Mi{\~n}arro.
\newblock A characterization of quasinormable {K}{\"o}the sequence spaces.
\newblock {\em Proc. Amer. Math Soc}, 123:1207--1212, 1995.

\bibitem{Moo10}
C.~G. Moorthy and C.~T. Ramasamy.
\newblock Characterizations of strictly singular and strictly discontinuous
  operators on locally convex spaces.
\newblock {\em Int. Journ. of Math. Analysis}, 4:1217--1224, 2010.

\bibitem{Mor00}
T.~J. Morrison.
\newblock {\em Functional Analysis: An Introduction to Banach Space Theory}.
\newblock Wiley \& Sons, 2000.

\bibitem{Ona91}
S.~{\"O}nal and M.~Yurdakul.
\newblock A note on strictly singular operators.
\newblock {\em Turk J Math}, 15(1):42--47, 1991.

\bibitem{Rya87}
R.~Ryan.
\newblock The {D}unford-{P}ettis property and projective tensor products.
\newblock {\em Bull. Polish. Acad. Sci. Math}, 35:785--792, 1987.

\bibitem{Rya02}
R.~Ryan.
\newblock {\em Introduction to the Tensor Products of Banach spaces}.
\newblock Springer-Verlag, 2002.

\bibitem{Ter86}
T.~Terzio{\~g}lu and M.~Yurdakul.
\newblock Restrictions of unbounded continuous linear operators on
  {F}r{\'e}chet spaces.
\newblock {\em Arch. Math}, 46:547--550, 1986.

\bibitem{Dul70}
D.~van Dulst.
\newblock Perturbation theory and strictly singular operators in locally convex
  spaces.
\newblock {\em Studia Math}, 38:341--372, 1970.

\bibitem{Dul71}
D.~van Dulst.
\newblock On strictly singular operators.
\newblock {\em Comp. Math}, 23:169--183, 1971.

\bibitem{Wro78}
V.~V. Wrobel.
\newblock Strikt singul{\"a}re {O}peratoren in lokalkonvexen {R}{\"a}umen.
\newblock {\em Math. Nachr}, 83:127--142, 1978.

\bibitem{Yur93}
M.~Yurdakul.
\newblock A remark on a paper of {J}. {P}rada.
\newblock {\em Arch Math}, 61:385--390, 1993.

\bibitem{Zah73}
V.~Zahariuta.
\newblock On the isomorphism of {C}artesian products of locally convex spaces.
\newblock {\em Studia Math.}, 46:201--221, 1973.

\end{thebibliography}

\end{document}